\documentclass[a4paper]{amsart}

\usepackage{amsmath,amsthm,amsfonts,amssymb}
\usepackage{graphicx}
\usepackage{latexsym}
\usepackage[all]{xy}
\usepackage{subfigure}
\usepackage[usenames,dvips]{color}  

\xyoption{matrix}
\xyoption{arrow}

\newtheorem{lemma}{Lemma}[section]

\newtheorem{theorem}[lemma]{Theorem}

\theoremstyle{definition}
\newtheorem*{remark}{Remark}
\newtheorem*{definition}{Definition}
\newtheorem{example}[lemma]{Example}
\newtheorem{problem}{Problem}

\begin{document}

\title[Cluster algebras and cluster categories]
{Geometric construction of 
cluster algebras and cluster categories}

\author[Baur]{Karin Baur}
\address{Department of Mathematics \\
ETH Z\"urich
R\"amistrasse 101 \\
8092 Z\"urich \\
Switzerland
}
\email{baur@math.ethz.ch}


\begin{abstract}
In this note we explain how to obtain cluster algebras  
from triangulations of (punctured) discs following the approach 
of~\cite{fst06}. 
Furthermore, we give a description of 
$m$-cluster categories via diagonals (arcs) in (punctured) polygons 
and of $m$-cluster categories via powers of translation quivers as given 
in joint work with R. Marsh (\cite{bm1}, \cite{bm2}). 
\end{abstract}

\maketitle

\tableofcontents

\section{Introduction}
This article is an expanded version of a talk presented at the 
Courant-Colloquium ``G\"ottingen trends in Mathematics'' in October 2007. 
It is a survey on two approaches to cluster algebras and 
($m$-)cluster categories via geometric constructions. 

Cluster algebras where introduced in 2001 by Fomin and Zelvinsky, 
cf.~\cite{fz}. They arose from the study of two related problems. 

\begin{problem}[Canonical basis]
Understand the {\em canonical basis} (Lusztig), 
or {\em crystal basis} (Kashiwara) of 
quantized enveloping algebras associated to a semisimple 
complex Lie algebra. It is expected that the positive part of the 
quantized enveloping algebra has a (quantum) cluster algebra structure, 
with the so-called cluster monomials forming part of the dual canonical basis. 

This picture motivated the definition of  {\em cluster 
variables}. 
\end{problem}

\begin{problem}[Total positivity]
An invertible matrix with real entries is called {\em totally positive} 
if all its minors are positive. This notion has been extended to all reductive 
groups by Lusztig~\cite{lu}. 
To check total positivity 
for an 
upper uni-triangular matrix, only a certain collection of the 
non-zero minors needs to be checked (disregarding the minors which are zero 
because of the uni-triangular from). The minimal sets of such all have the same cardinality. 
When one of them is removed, it can often be replaced by a unique 
alternative minor. 
The two minors are connected through a certain relation. 

This exchange ({\em mutation for minors}) motivated 
the definition of cluster mutation.
\end{problem}

The subject of cluster algebra is a very young and dynamic one. In 
the past few years, connections to various other fields arose. 
We briefly mention a few of them here. 

\begin{itemize}
\item
Poisson geometry (integrable systems), Teichm\"uller spaces (local coordinate 
systems), cf. Gekhtman-Shapiro-Vainshtein~\cite{gsv1,gsv2} and 
Fock-Goncharov~\cite{fg}; 
\item $Y$-systems in thermodynamic Bethe Ansatz 
(families of rational functions defined by recurrences which were 
introduced by Zamolodchikov~\cite{za}). 
Cf.~\cite{fz}; 
\item Stasheff polytopes, associahedra, Chapoton-Fomin-Zelevinsky~\cite{cfz};
\item ad-nilpotent ideals of Borel subalgebras in Lie algebras, Panyushev~\cite{pa};
\item Preprojective algebra models, Geiss-Leclerc-Schr\"oer,~\cite{gls1}, \cite{gls2};
\item Representation theory, tilting theory, etc., Cf. e.g.~\cite{bmrrt}.
\end{itemize}

In this article, we will first recall triangulations of surfaces with 
marked points and associate certain integral valued matrices to them. 
Then we will give a brief introduction to cluster algebras (Section~\ref{s:clusteralg}). 
In Section~\ref{s:triang-clu_a} we show how to associate cluster algebras 
to triangulations of (punctured) discs. 
Then we explain what cluster categories and $m$-cluster categories are 
(Section~\ref{s:clustercat}) 
and give a combinatorial model to describe $m$-cluster categories via 
arcs in a polygon in Section~\ref{s:quivers-categories}, cf. 
Theorems~\ref{thm:A},~\ref{thm:B}
as given in our 
joint work with R. Marsh (\cite{bm1},~\cite{bm2}). 
In addition, we 
obtain a descriptions of the $m$-cluster categories 
using the notion of the power of a translation quiver (Theorem~\ref{thm:power}).  
At the end we describe connections to other work, pose several questions 
and show new directions in 
this young and dynamic field (Section~\ref{s:connection}). 

%
\section{Triangulated surfaces}\label{s:triang}
%
%
In this section we recall triangulation of surfaces following the approach 
of Fomin, Shapiro and Thurston~\cite{fst06}. 
Let $S$ be a connected oriented Riemann surface with 
boundary. Fix a finite set $M$ of {\it marked points} on 
$S$. Marked points in the interior of $S$ are called 
{{\it punctures}}. 

We consider triangulations of $S$ whose vertices are 
at the marked points in $M$ and whose edges 
are pairwise non-intersecting curves, 
so-called {\it arcs}
connecting marked points. The most important example for us is the 
case where $S$ is a disc with marked points on the boundary 
and with at most one puncture. We will later restrict to that case but 
for the moment we explain the general picture. 

It is convenient to exclude cases where there are no such 
triangulations (or only one such). We always assume that 
$M$ is non-empty and that each boundary component has at least one 
marked point. 
And we disallow the cases $(S,M)$ with one 
boundary component, $|M|=1$ with $\le 1$ puncture and $|M|\in\{2,3\}$ 
with no puncture. 

In case $S$ is a (punctured) disc we will also call it a {\em (punctured) polygon}. 
E.g. if $(S,M)$ has three marked points on the boundary and a puncture, we 
will say that $S$ is a once-punctured triangle. 

\begin{figure}
\begin{center}
\subfigure[Once-punctured triangle]
	{\includegraphics[scale=.5]{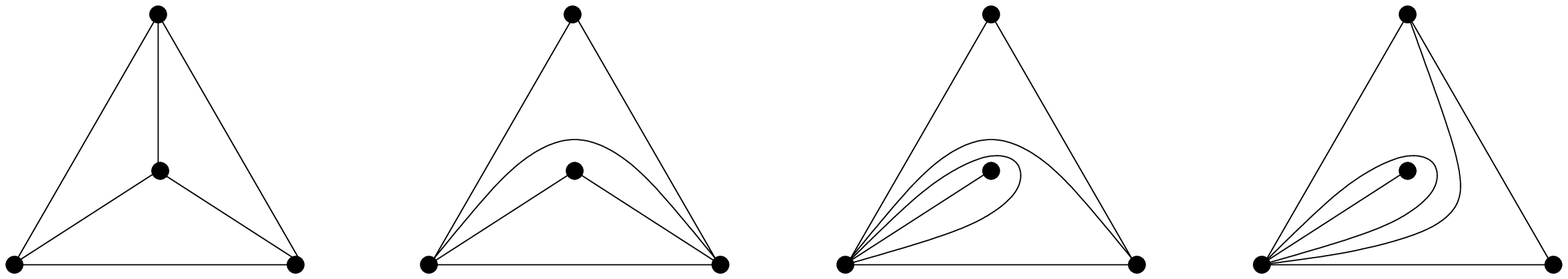}\label{fig:triangle}}
	\\
\subfigure[Annulus]
	{\includegraphics[scale=.5]{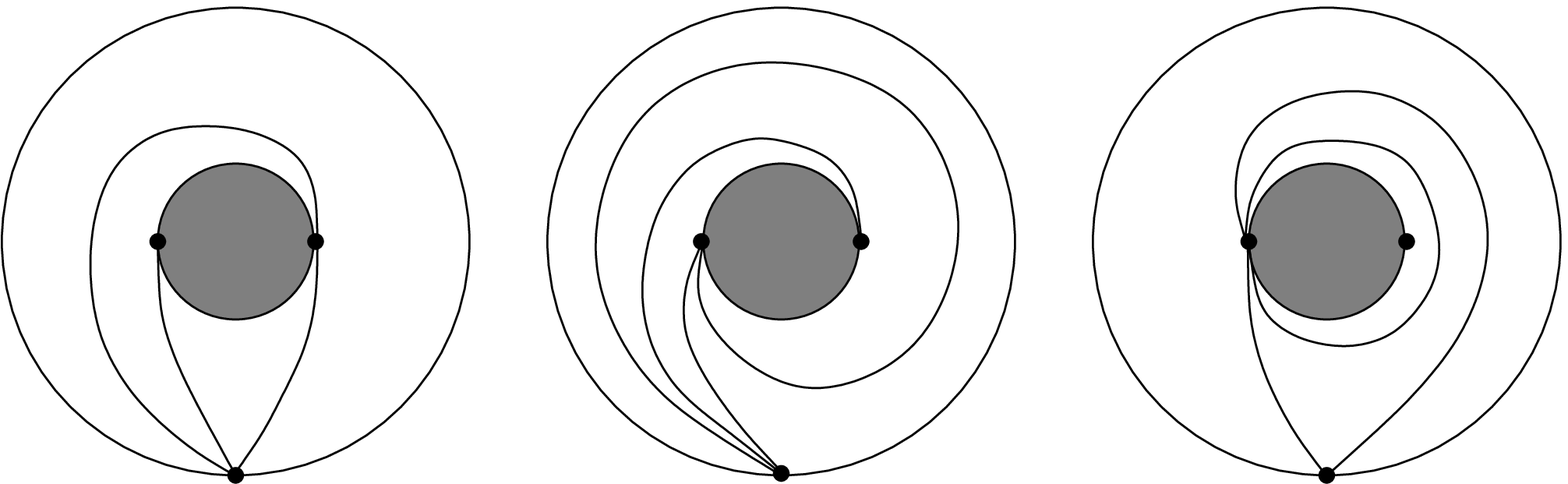}\label{fig:annulus}}
\end{center}
\caption{Examples of triangulations}\label{fig:triangulations}
\end{figure}

Note that the pair $(S,M)$ is defined (up to homeomorphism) 
by the genus of $S$, by the numbers of boundary components, 
of marked points on each boundary component 
and of punctures. 
Two examples of such triangulations are given in Figure~\ref{fig:triangulations}. 

\begin{definition}
A curve in $S$ (up to isotopy relative $M$) 
is an {\it arc $\gamma$ in $(S,M)$} if 
\begin{tabular}{ll}
(i) & the endpoints of $\gamma$ are marked points in $M$; \\
(ii) & $\gamma$ does not intersect itself (but its endpoints might coincide);\\
(iii) & relative interior of $\gamma$ is disjoint from $M$ and 
from the boundary of $S$; \\
(iv) & $\gamma$ does not cut out an unpunctured monogon or digon.
\end{tabular}
\end{definition}

The set of all arcs in $(S,M)$ is usually infinite as we can already see in 
the case of the annulus of Figure~\ref{fig:triangulations}(b). 
One can show that it is finite if and only 
if $(S,M)$ is a disk with at most one puncture, i.e. if $(S,M)$ is the object of 
our interest. 

Two arcs are said to be {\it compatible} 
if they do not intersect in the interior of $S$. An 
{\it ideal triangulation} is a maximal 
collection $T$ of pairwise compatible arcs. The arcs of $T$ cut $S$ 
into the so-called {\it ideal triangles}. 
These triangles may be self-folded, e.g. along the horizontal arc in the 
picture below: 

\begin{center}
\includegraphics[scale=.6]{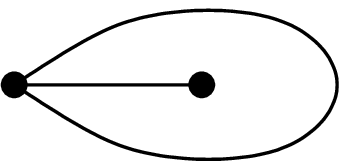}
\end{center}

An easy count shows that the once-punctured triangle has 
ten ideal triangulations, the four of figure~\ref{fig:triangulations}, with 
the rotations of the last three (by $120^{\circ}$ and $240^{\circ}$). 

In fact we can say more: the number of arcs in an ideal triangulation 
is an invariant of $(S,M)$, we call it the 
{\it rank} of $(S,M)$. There is a formula for it, cf.~\cite{fg2}: 
if $g$ is the genus of $S$, 
$b$ the number of boundary components, 
$p$ the number of punctures, $c$ the number of marked points on the boundary, 
then the rank of $(S,M)$ is
\[
6g+3b+3p+c-6
\]
The rank of the once punctured triangle of Figure~\ref{fig:triangulations}(a) 
is thus three as expected. 

For small rank,~\cite[Example 2.12]{fst06} gives a list of all 
possible choices of $(S,M)$. The word ``type'' appearing in the 
list refers 
to the Dynkin type of to the corresponding cluster algebra as will  
be explained later: 

\begin{tabular}{ll}
Rank 1 & unpunctured square (type ${\rm A}_1$) \\
Rank 2 & 
unpunctured pentagon (type ${\rm A}_2$) \\
 & once-punctured digon (type ${\rm A}_1\times {\rm A}_1$) \\
 & annulus with one marked point on each boundary component \\
Rank 3 & 
unpunctured hexagon (type ${\rm A}_3$) \\
& once-punctured triangle (type ${\rm A}_3={\rm D}_3$) \\
& annulus with one marked point on one boundary component, two 
on the other \\
& once-punctured torus.
\end{tabular}

If $T$ is an ideal triangulation of $(S,M)$ and $p$ an arc of $T$ 
as in the picture below, 
we can replace $p$ by an arc $p'$ through a so-called flip or Whitehead 
move: 

\begin{center}
\includegraphics[scale=0.6]{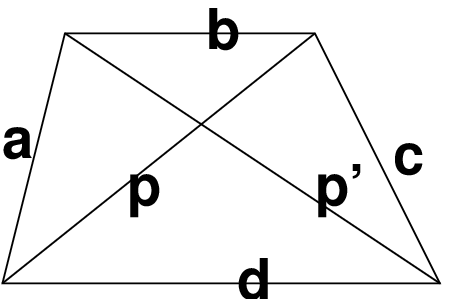}
\end{center}

Here we allow that some of the sides $\{a,b,c,d\}$ coincide. 
A consequence of a result of Hatcher (\cite{ha}) 
is that for any two ideal triangulations 
$T$ and $T'$ there exists a sequence of flips 
leading from $T$ to $T'$.

We next want to associate a matrix to an ideal triangulation of 
$(S,M)$. This works as follows. 
Let $T$ be an ideal triangulation of $(S,M)$, label the arcs of $T$ 
by $1$, $2$, $\dots,n$. Then define 
$B(T)$ to be the following $n\times n$-square matrix 
\[
B(T)=\sum_{\Delta} B^{\Delta}
\]
where the $n\times n$-matrices $B^{\Delta}$ are defined for 
each triangle $\triangle$ of $T$ by 
\[
b_{ij}^{\Delta} \ = 
\ \left\{
\begin{array}{rl}
1     & \mbox{if $\Delta$ has sides 
{\it i} and {\it j} where {\it j} is a clockwise 
neighbour of {\it i};} \\
 -1  & \mbox{if $\Delta$ has sides {\it i} and {\it j} where 
{\it i} is a clockwise 
neighbour of {\it j};} \\
0    & \mbox{otherwise.}
\end{array}
\right.
\]
The matrix $B(T)$ is skew-symmetric with entries $0,\pm 1, \pm 2$. 

\begin{figure}
\begin{center}
\includegraphics[scale=0.6]{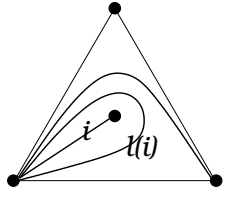}
\end{center}
\caption{Enclosing loop $l(i)$ of the arc $i$}\label{fig:bloop}
\end{figure}

\begin{remark}
In order to simplify the definition of $b_{ij}^{\Delta}$ we have cheated a little bit. 
Whenever the triangle $\Delta$ is self-folded along an arc $i$, then in the 
right hand side of the definition of the entry  $b_{ij}^{\Delta}$, the arc $i$ 
has to be replaced by its enclosing loop $l(i)$, 
cf. Figure~\ref{fig:bloop}. 
\end{remark}

\begin{example}
(1) We compute $B(T)$ for the triangulated 
punctured triangle. 
\begin{center}
\includegraphics[scale=0.4]{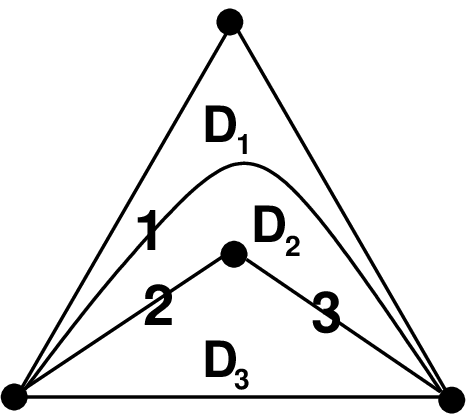}
\end{center}
It is: 
\[
\underbrace{B^{D_1}}_{=\ 0}+B^{D_2} + B^{D_3} \ = \  
\begin{pmatrix}
0&1&-1\\
-1&0&1\\
1&-1&0
\end{pmatrix}
+ 
\begin{pmatrix}
0&0&0\\
0&0&-1\\
0&1&0
\end{pmatrix}
=
\begin{pmatrix}
0&1&-1\\
-1&0&0\\
1&0&0
\end{pmatrix}
\]

(2)
Take an annulus with one marked point on each 
boundary and the triangulation $T$ as in the picture. Then $B(T)$ is 
\[
B^{D_1} + B^{D_2} \ = \  
\begin{pmatrix}
0&1\\
-1&0&
\end{pmatrix}
+ 
\begin{pmatrix}
0&1\\
-1&0
\end{pmatrix}
=
\begin{pmatrix}
0&2\\
-2&0
\end{pmatrix}
\]
\begin{center}
\includegraphics[scale=.4]{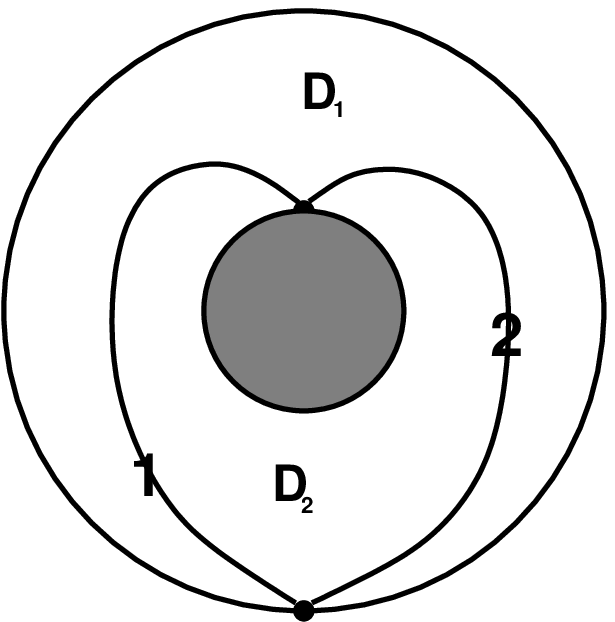}
\end{center}
\end{example}

%
\section{Cluster algebras}\label{s:clusteralg}
%
In this section we present a very short introduction to cluster algebras, 
following Fomin-Zelevinsky~\cite{fz}. 
A cluster algebra $\mathcal A=\mathcal{A}(\underline{x},B)$ is a 
subring of $\mathbb{F}=\mathbb{Q}(u_1,\dots,u_m)$, associated 
to a {\it seed} $(\underline{x},B)$ defined in the following way. 
\begin{itemize}
\item[(i)] A seed is a pair $(\underline{x},B)$ consisting of a 
{\em cluster} $\underline{x}=(x_1,\dots,x_m)$ where $\underline{x}$ 
is a free generating set of $\mathbb{F}$ over $\mathbb{Q}$ 
and 
$B=(b_{xy})_{xy}$ is a sign-symmetric $m\times m$ matrix with integer 
coefficients, 
i.e. $b_{xy}\in\mathbb{Z}$ for all $1\le x,y\le m$ and if 
$b_{xy}>0$ then $b_{yx}<0$. 
\item[(ii)]  A seed $(\underline{x},B)$ can be mutated to another seed 
$(\underline{x'},B')$: {\em mutation at $z\in\underline x$} is 
the map $\mu_z$: 
$(\underline{x},B)$ $\mapsto\ (\underline{x'},B')$: 
$\underline{x'}= \underline{x} -\,z\ \cup z'$ where $z'$ is defined via 
the {\em exchange relation} 
\[
zz' \ = \ \prod_{\stackrel{x\in\underline{x}\ \ }{b_{xz}>0}}x^{b_{xz}} \ +
\ \prod_{\stackrel{x\in\underline{x}\ \ }{b_{xz}<0}}x^{-b_{xz}}
\]
and $B'$ is defined similarly via {\em matrix mutation}: 
\[
b_{xy}' \ = \left\{
    \begin{array}{ll}
     -b_{xy} & \mbox{if $x=z$ or $y=z$}\\
      b_{xy}+\frac{1}{2}(|b_{xz}|\cdot b_{zy}+b_{xz}\cdot |b_{zy}|) & \mbox{otherwise.}
     \end{array}
\right.
\]
\end{itemize}
Note: $\mu$ is involutive, i.e. $\mu_{z'}(\mu_z((\underline{x},B)))=(\underline{x},B)$.

\noindent
Two seeds $(\underline{x},B)$ and $(\underline{x'},B')$
are said to be {\em mutation-equivalent} if one can be 
obtained from the other through a sequence of mutations. \\
The {\it cluster variables} are defined to be 
the union of 
all clusters of a mutation-equivalence 
class (of a given seed). These appear in overlapping sets. 
Finally, the corresponding 
{\it cluster algebra} $\mathcal A=\mathcal{A}(\underline{x},B)$ 
is the subring of $\mathbb{F}$ generated by all the cluster variables. 
(Here we are defining cluster algebras with trivial coefficients.)
A cluster algebra is said to be {\em of finite type} if there exists only a 
finite number of cluster variables. 

One can show that up to isomorphism of cluster algebras 
$\mathcal A(\underline{x},B)$ 
does not depend on the initial choice of a free generating set 
$\underline{x}$.

\begin{example}\label{ex:mutation}
(Type ${\rm A}_2$). 
We start with the pair $\underline{x}=(x_1,x_2)$, 
B$=\begin{pmatrix}
0&1\\
-1&0
\end{pmatrix}$. 
In a first step we mutate $x_1$. 
from $x_1x_1' = 1 + x_2$ we obtain 
$x_1' = \frac{1+x_2}{x_1}$. The next mutation is at $x_2$ (mutation at $x_1'$ 
would lead us back to $x_1$), we have 
$x_2'=\frac{x_1+1}{x_2}$. 
And then 
$x_1''= \frac{x_1+x_2+1}{x_1x_2}$;   
$x_2''=x_1$, $x_1'''=x_2$. \\
In particular, we obtain five cluster variables in this example. 
\end{example}

Some of the main results on cluster algebras are summarized here: 
\begin{itemize}
\item
{\it Laurent phenomenon:} 
${\mathcal A}(\underline{x},B)$ sits inside $\mathbb{Z}[x_1^{\pm},\dots,x_m^{\pm}]$, 
i.e. every element of the cluster algebra is an integer Laurent polynomial 
in the variables of $\underline{x}$ (cf.~\cite{fz-Laurent});
\item
Classification of finite type cluster algebras by roots systems,~\cite{fz-finite} 
(cluster algebras of finite type can be classified by Dynkin diagrams); 
\item
Realizations of algebras of regular functions on double 
Bruhat cells in terms of cluster algebras (\cite{bfz}). 
\end{itemize}

Examples of cluster algebras are: 
Coordinate rings of 
${\rm SL}_2$, ${\rm SL}_3$ (\cite{fz-notes}); 
Pl\"ucker coordinates on ${\rm Gr}_{2,n+3}$ (\cite{sco}, \cite{gsv1}). 

%
\subsection*{Cluster algebras and quivers}
%
%
We will now explain how to associate a quiver to a seed of a 
cluster algebra. 

Recall that a {\it quiver $\Gamma=(\Gamma_0,\Gamma_1)$} 
is an oriented graph with 
vertices $\Gamma_0$ and arrows $\Gamma_1$ between them.  
E.g. 
\[
1\stackrel{\alpha} {\longrightarrow}2\stackrel{\beta}{\longrightarrow} 3
\]
with $\Gamma_0=\{1,2,3\}$ and $\Gamma_1=\{\alpha,\beta\}$; 

Any skew-symmetric $m\times m$-matrix $B$ determines a quiver 
$\Gamma(B)$ 
with $m$ vertices. One labels the columns of $B$ by 
$\{1,2,\dots,m\}$ and sets $\Gamma_0=\{1,2,\dots,m\}$. 
Then one draws $b_{xy}$ arrows from $x$ to $y$ if $b_{xy}>0$ 
(for $x$, $y\,\in\ \Gamma_0$).

Such a quiver has no loops and for any two vertices $i\neq j$ of 
$\Gamma(B)$, there are only arrows in one direction between them. 

So in particular, if the matrix $B$ of a seed $(\underline{x}, B)$ is skew-symmetric, 
it determines a quiver in this way. 

\begin{example}
The matrix B$=\begin{pmatrix}
 0 & 1 \\
 -1 & 0
 \end{pmatrix}$
from Example~\ref{ex:mutation} above gives the quiver: 
\[
1\longrightarrow 2
\]
\end{example}

Clearly, this process is reversible: a quiver whose arrows only go in one direction 
between any given pair $i\neq j$ of vertices and without loops gives rise to 
a skew-symmetric matrix which we will denote by $B(\Gamma)$.

%
\section{From triangulations to cluster algebras}\label{s:triang-clu_a}
%
From now on we assume that $(S,M)$ is a disc with at most one puncture. 
We want to show how a triangulation $T$ of $(S,M)$ determines a cluster algebra. 
Label the arcs of $T$ by $1,2,\dots,n$. 

%
Then we define a cluster $\underline{x}_T=(x_1,\dots,x_n)$ 
by sending $i\mapsto x_i$ and choose as a matrix the 
the skew-symmetric matrix associated $B(T)$ associated to $T$. 
in Section~\ref{s:triang}. 
This clearly produces a seed 
$(\underline{x}_T,B(T))$. 

Thus to the triangulation $T$ of the disc $(S,M)$ we have associated the seed 
$(\underline{x}_T,B(T))$ and hence obtain a cluster 
algbra $\mathcal{A}=\mathcal{A}(\underline{x}_T,B(T))$.

\begin{example}
Consider an unpunctured pentagon as below. In the triangulation, we label the 
arcs $1$ and $2$. They form a triangle $D_1$ 
together with a boundary arc 
and $2$ is the clockwise neighbour of $1$. 

\begin{center}
\includegraphics[scale=.4]{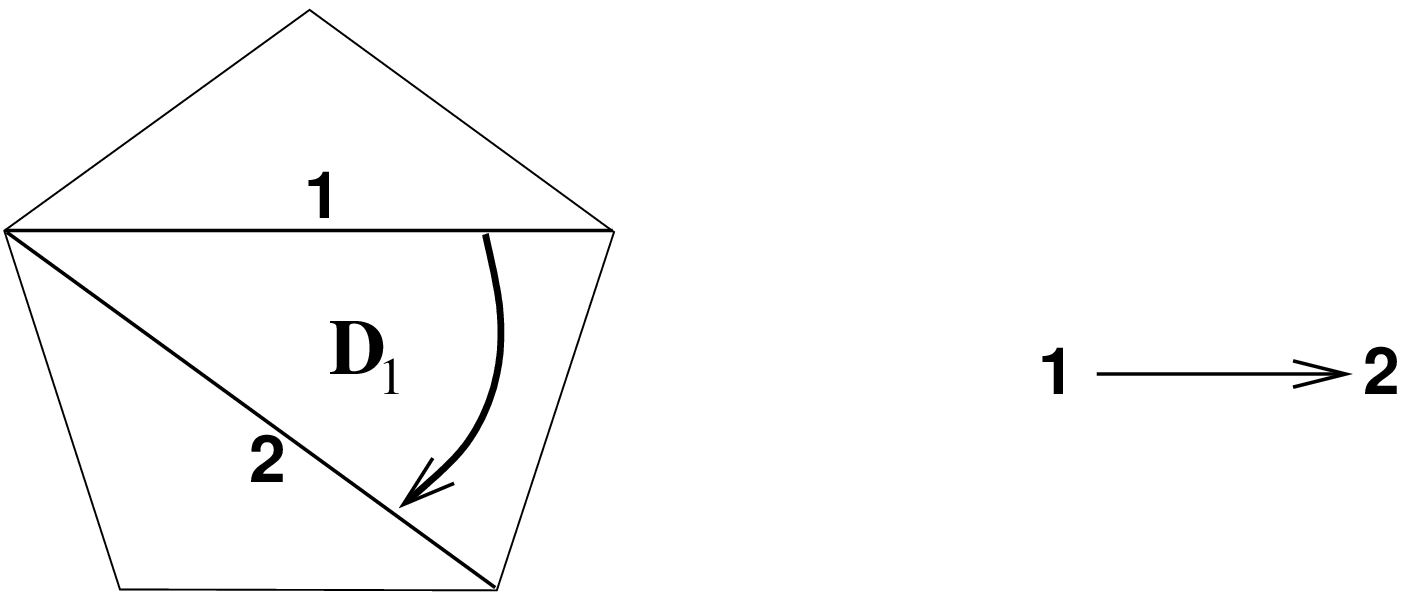} 
\end{center}
Then the seed we obtain is $((x_1,x_2),B(T))$ with 
$B(T)=\begin{pmatrix}
0 & 1 \\ -1 & 0
\end{pmatrix}$ 
as in Example~\ref{ex:mutation} above. 
\end{example}

%
\section{Cluster categories}\label{s:clustercat}
%

Cluster categories are certain quotients of derived categories 
of module categories. 
They were 
introduced in 2005 by 
Buan-Marsh-Reineke-Reiten-Todorov (\cite{bmrrt}). 

Independently, Caldero-Chapoton-Schiffler have introduced 
the cluster categories (in type $A_n$) in 2005 (\cite{ccs}) 
using a graphical description. Later, Schiffler extended this 
to type $D_n$ in (\cite{schi}). 

The aim behind the definition of cluster categories was to model 
cluster algebras using the representation theory of quivers. This 
was motivated by the observation that the cluster variables of a cluster 
algebra of finite type are parametrized by the almost positive roots of 
the corresponding root system. 

Cluster categories 
have led to new development in the theory 
of the (dual of the) canonical bases, they provide insight into cluster 
algebras. They have 
also developed into a field of their own. E.g. they have 
led to the definition of cluster-tilting theory.

Let us describe the construction of cluster categories, following~\cite{bmrrt}. 

We start with a quiver $Q$ whose underlying graph is a simply-laced Dynkin 
diagram (i.e. of type ADE). Denote by 
$D^b(kQ)$ the bounded derived category of finite dimensional 
$kQ$-modules (we assume that the field $k$ is algebraically closed). 
Note that the shape of the quiver of $D^b(kQ)$ is $Q\times \mathbb{Z}$ 
with certain connecting arrows. By {\em quiver of $D^b(kQ)$} we mean 
the Auslander-Reiten quiver of $D^b(kQ)$, i.e. the quiver whose vertices 
are the isomorphism 
classes of indecomposable modules and whose arrows come from irreducible 
maps between them. 

This quiver has two well-known graph automorphisms: 
$\tau$ (``Auslander-Reiten translate'') which sends each vertex to its 
neighbour to the left. And 
$[1]$ (the ``shift'') which sends a vertex in a copy of the module category 
of $kQ$ to the corresponding vertex in the next copy of the module category. 




The {\em cluster category}, $\mathcal{C}$, is now defined as the orbit 
category of $D^b(kQ)$ under a canonical automorphism: 
\[
\mathcal{C} :=\mathcal{C}(Q) :=D^b(kQ)/\tau^{-1}\circ[1] 
\]
One can show that this is independent of the chosen orientation of $Q$. 
More generally, Keller (\cite{ke}) has introduced the $m$-cluster category, 
{$\mathcal{C}^m$} as follows: 
\[
\mathcal{C}^m :=\mathcal{C}^m(Q):=D^b(kQ)/\tau^{-1}\circ[m]
\]
Keller has shown in \cite{ke} that $\mathcal{C}^m$ is triangulated 
and a Calabi-Yau category of dimension $m+1$. 
Furthermore, $\mathcal{C}^m$ is Krull-Schmidt (\cite{bmrrt}). 
The $m$-cluster category has attracted a lot of interest over the last few 
years. In particular, it has been studied by Keller-Reiten, 
Thomas, Wralsen, Zhu, B-Marsh, 
Assem-Br\"ustle-Schiffler-Todorov, Amiot, Wralsen, etc.

Our goal for this note is to describe $\mathcal{C}^m$ using diagonals 
of a polygon 
(type ${\rm A}_n$) and arcs in a punctured polygon (type ${\rm D}_n$).

%
\section{From arcs via quivers to cluster categories}\label{s:quivers-categories}
%
%
Let us first recall the notion of a stable translation quiver due to 
Riedtmann~\cite{rie}. 

\begin{definition}
A {\it stable translation quiver} is a 
pair $(\Gamma,\tau)$ where 
$\Gamma=(\Gamma_0,\Gamma_1)$ is a quiver (locally finite, without loops)
and $\tau:\Gamma_0\to\Gamma_0$ is a bijective map
such that the number of arrows from $x$ to $y$ equals the number 
of arrows from $\tau y$ to $x$ 
%
for all $x,\ y\in\Gamma_0$. 
The map $\tau$ is called the {\it translation} 
of $(\Gamma,\tau)$. 
\end{definition}

Now we are ready to define a quiver $\Gamma$ from a hexagon (see 
figure below) as follows: 
\begin{center}
\includegraphics[scale=.7]{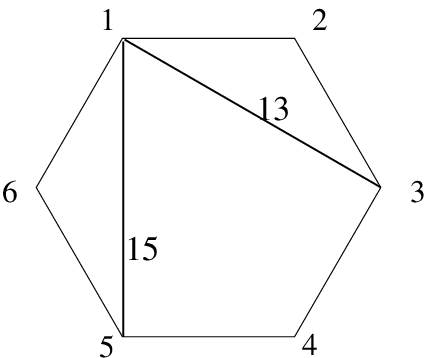}
\end{center}
\noindent
$\Gamma_0$: 
The vertices are the diagonals $(ij)$ of the hexagon 
($1\le i<j-1\le 7$). \\
$\Gamma_1$: 
The arrows are of the form 
$(ij)\to (i,j+1)$, $(ij)\to (i+1,j)$
provided the target is also a diagonal in the hexagon ($i,j$ $\in\mathbb{Z}_6$). \\
Set $\tau$: $(ij)\mapsto (i-1,j-1)$ to be anti-clockwise rotation 
about the center. \\

\noindent
The quiver obtained this way from the hexagon is the following: 
\[
\small
\xymatrix@-4mm{
 & & 15\ar[rd]\ar@{--}[rr] & & 26\ar[rd]\ar@{--}[rr] & & 13\ar[rd]\\
 & 14\ar[ru]\ar[rd]\ar@{--}[rr] & & 25\ar[rd]\ar[ru]\ar@{--}[rr]
 & &  36\ar[ru]\ar[rd]\ar@{--}[rr] & & 14\ar[rd]\\
13\ar@{--}[rr]\ar[ru] & & 24\ar@{--}[rr]\ar[ru] & & 35\ar[ru]\ar@{--}[rr]
 & & 46\ar[ru]\ar@{--}[rr] & & 15}
\]
It clearly is an example of a stable translation quiver. 

Note that such a quiver can be defined for any polygon. Denote 
the quiver arising in that way by $\Gamma(n,1)$ if $n+2$ is the number 
of vertices of the polygon. (The use of $n$ instead of $n+2$ in the notation of 
the quiver $\Gamma(n,1)$ and the 
extra entry $1$ are used to make this compatible with the more general 
setting involving $m$-diagonals described below). 
Caldero, Chapton 
and Schiffler have shown that the cluster category can be obtained via 
diagonals in a polygon: 

\begin{theorem}[\cite{ccs}]
The quiver of the cluster category $\mathcal{C}=\mathcal{C}(A_{n-1})$ is 
isomorphic to the quiver $\Gamma(n,1)$ obtained from an $(n+2)$-gon.
\end{theorem}
As before in the case of the bounded derived category, 
the {\em quiver of $\mathcal{C}$} is an abbreviation for the 
Auslander-Reiten quiver of $\mathcal{C}$. It has as vertices the indecomposable 
objects of $\mathcal{C}$, and as arrows are the 
irreducible maps between them. 

To be able to model $m$-cluster categories we now generalize the notion 
of diagonal and introduce the so-called $m$-diagonals. We start with 
a polygon $\Pi$ with $nm+2$ vertices ($n,m\in\mathbb{N}$), labeled by 
$1,2,\dots,nm+2$. 

\begin{definition}
An {\it $m$-diagonal} is a diagonal 
$(ij)$ dividing $\Pi$ into
an $mj+2$-gon and an $m(n-j)+2$-gon ($1\le j\le\lceil\frac{n-1}{2}\rceil$). \\
\end{definition}

\begin{example} To illustrate this, let $\Pi$ be an octagon, $n=3$, $m=2$. 
In that case, $1\le j\le 1$, so any $2$-diagonal has to divide $\Pi$ into 
a quadrilateral and a hexagon. 
\begin{center}
\includegraphics[scale=.8]{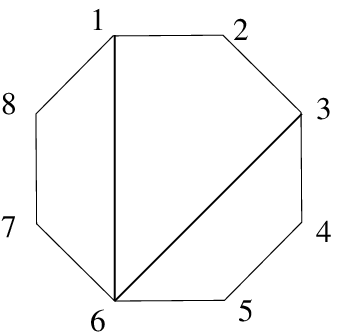}
\end{center}
Observe that each maximal set of non-crossing $2$-diagonals 
contains two elements. They are 
$\{(16),(36)\}$, $\{(16),(25)\}$, $\{(16),(14)\}$ 
and rotated version of these. 
\end{example}
Recall that the number of arcs in a triangulation (see Section~\ref{s:triang}) 
is an invariant of a disc $(S,M)$, called the rank of $(S,M)$. 
In the same way, 
the maximal number of non-crossing $m$-diagonals is an invariant 
of the polygon. It is equal to $n-1$ (for the $nm+2$-gon $\Pi$). 

Using $m$-diagonals we can now define a translation quiver 
$\Gamma(n,m)=(\Gamma,\tau_m)$:

{\it $\Gamma_1$: $(ij)\to (ij')$\ \ } if $(ij)$, 
$B_{jj'}$ and $(ij')$ span an $m+2$-gon 
($B_{jj'}$ is the boundary $j$ to $j'$, going clockwise).

\noindent
$\Gamma_0$: 
The vertices are the $m$-diagonals $(ij)$ in $\Pi$ (with $1\le i<j-m$).\\
$\Gamma_1$: 
The arrows are of the form 
$(ij)\to (i,j+m)$, $(ij)\to (i+m,j)$ provided the target is still inside 
the polygon. In other words: $(ij)$, $(i,j+m)$ and the boundary 
arc $j$ to $j+m$ (resp. $(ij)$, $(i+m,j)$ and the boundary arc 
from $i$ to $i+m$) form an $m+2$-gon as in the picture: 
\begin{center}
\includegraphics[scale=.4]{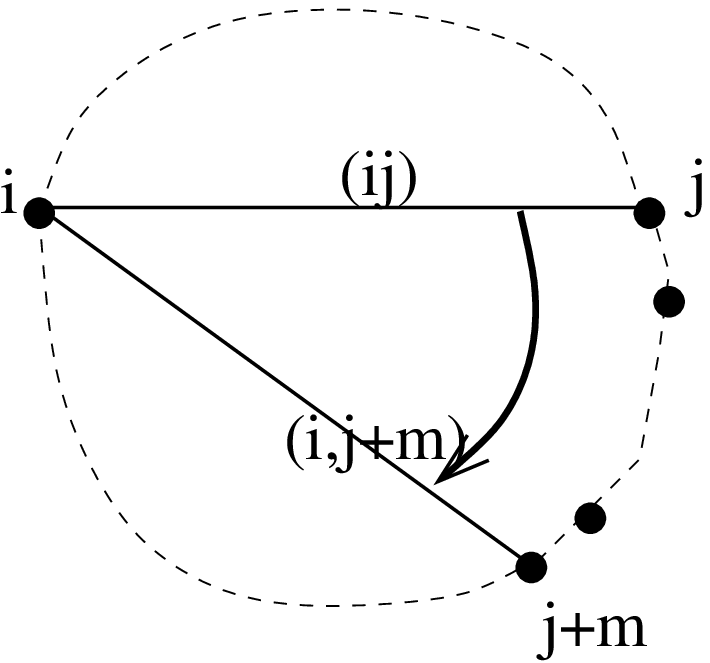}
\end{center}
Furthermore, let $\tau_m$ be anti-clockwise rotation (about center) through 
the angle $m\frac{2\pi}{nm+2}$. 

\begin{remark}
It is clear that 
$\Gamma(n,m)$ is a stable translation quiver. 
In case $m=1$, we recover the usual diagonals.
\end{remark}

The quiver $\Gamma(3,2)$ for the octagon from the previous example 
is thus: 
\[
\small
\xymatrix@-6mm{
 & 16\ar[rdd]\ar@{--}[rr] && 38\ar[rdd]\ar@{--}[rr] && 25\ar[rdd]\ar@{--}[rr]
 && 47\ar[rdd]\ar@{--}[rr] && 16\\
 & & \\
14\ar@{--}[rr]\ar[ruu] & & 36\ar[ruu]\ar@{--}[rr] & & 58\ar[ruu]\ar@{--}[rr]
 & &  27\ar[ruu]\ar@{--}[rr] & & 14\ar[ruu]}
\]

Then one can show that the $A$-type $m$-cluster category can be obtained using 
$m$-diagonals in a polygon: 

\begin{theorem}[\cite{bm1}]\label{thm:A}
The quiver of the $m$-cluster category $\mathcal{C}^m={\mathcal C}^m(A_{n-1})$ 
is isomorphic to the quiver $(\Gamma(n,m),\tau_m)$ obtained from $m$-diagonals 
in an $nm+2$-gon. 
\end{theorem}

The proof of our result uses Happel's description of the Auslander-Reiten-quiver of 
$D^b(kQ)$ where $Q$ is of Dynkin type $A_{n-1}$ 
and combinatorial analysis of $\Gamma(n,m)$. For details we refer 
to~\cite[Section 5]{bm1}. 

%
\subsection*{The description in type $D$}
%

We have a similar description of the $m$-cluster categories 
of $D$-type. Instead of working with a polygon (or unpunctured disc) 
we now have to use a punctured polygon. Let $\Pi$ be a punctured 
$nm-m+1$-gon. Instead of using the term diagonal, we now 
speak of {\em arcs} in $\Pi$. An arc going from $i$ to $j$, homotopic 
equivalent to the boundary $B_{ij}$ from $i$ to $j$ (going clockwise) is 
denoted by $(ij)$. By $(ii)$ we denote an arc homotopic equivalent to the 
boundary $B_{ii}$ with endpoints in $i$. And $(i0)$ is an arc homotopic 
equivalent to the arc between $i$ and the puncture $0$. 
We will say that an $n$-gon is {\em degenerate} if 
it has $n$ sides and $n-1$ vertices. 

For details and examples we refer to~\cite[Section 3]{bm2}.

\begin{definition}
An {\em $m$-arc} of $\Pi$ is an arc $(ij)$ such that 

(i) $(ij)$ and $B_{ij}$ (the boundary from $i$ to $j$, going clockwise) form 
an $km+2$-gon for some $k$,

(ii) $(ij)$ and $B_{ji}$ (the boundary from $j$ to $i$, going clockwise) form 
an $lm+2$-gon for some $l$. 

Furthermore, $(ii)$ and $(i0)$ are called {\em $m$-arcs} if $(ii)$ and $B_{ii}$ form 
a degenerate $km+2$-gon for some $k$. 
\end{definition}

Then we can define a translation quiver 
$\Gamma=\Gamma_{\odot}(n,m)$ 
as follows: 

\noindent
$\Gamma_0$: The vertices are the $m$-arcs of $\Pi$ \\
$\Gamma_1$: The arrows are the so-called $m$-moves between 
vertices: \\

We say that $(ij)\to (ik)$ is an {\em $m$-move} if $(ij)$, $B_{jk}$  
and $(ik)$ span a (degenerate) $m+2$-gon. In the figure below 
there are two examples of $2$-moves. 
\begin{center}
\includegraphics[scale=.4]{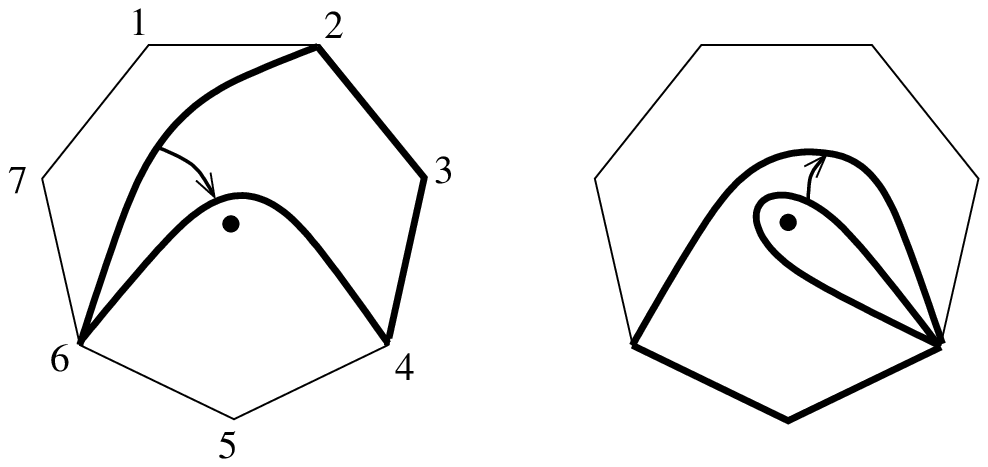}
\end{center}
$\tau_m$: rotation anti-clockwise (about center), through an angle of 
$\frac{m}{nm-m+1}$. 

Clearly, the pair $(\Gamma_{\odot}(n,m),\tau_m)$ is also a stable translation quiver. 
We can now formulate the statement. 

\begin{theorem}[Theorem 3.6 of \cite{bm2}]\label{thm:B}
The quiver of the $m$-cluster category $\mathcal{C}^m={\mathcal C}^m(D_n)$ 
is isomorphic to the quiver $(\Gamma_{\odot}(n,m),\tau_m)$ obtained 
from $m$-arcs in an $nm-m+1$-gon. 
\end{theorem}

%
\subsection*{The $m$-th power of a translation quiver}
%
%
We will now describe another way to obtain $m$-cluster categories directly 
from the diagonals or arcs in a (punctured) polygon. 
Let $(\Gamma,\tau)$ a translation quiver as before. 
Then we define the {\em $m$-th power of $\Gamma$}, $\Gamma^m$, 
to be the quiver whose vertices are the vertices of $\Gamma$ 
(i.e. $\Gamma^m_0=\Gamma_0$). 
The arrows in $\Gamma^m$ are 
paths of length $m$, going in an unique direction. 
(To be precise, we ask that such a path is {\em sectional}, i.e. 
that in a path $x_0\to x_1\to\dots\to x_{m-1}\to x_m$ of length $m$ 
we have $\tau x_{i+1}\neq x_{i-1}$ whenever $\tau x_{i+1}$ 
is defined.) 
And the translation $\tau^m$ of $\Gamma^m$ is obtained by repeating 
the original translation $m$ times. 

\begin{definition}
The quiver $(\Gamma^m,\tau^m)$ as defined above is called the 
{\em $m$-th power of $(\Gamma,\tau)$}. 
\end{definition}

With this we are ready to formulate the result: 

\begin{theorem}[\cite{bm1}]\label{thm:power}
$\Gamma(n,m)$ is a connected component of
$(\Gamma(nm,1))^m=(\Gamma(\text{cluster category}))^m$. 
\end{theorem}

\begin{remark}
Observe that $(\Gamma^m,\tau^m)$ is again a stable translation quiver. 
However, even if $(\Gamma,\tau)$ is connected, the $m$-th power is 
not connected in general! 
\end{remark}

\begin{example}
To illustrate this consider the quiver $\Gamma(6,1)$ of an octagon. 
\[
\small
\xymatrix@-5mm{
 & & & & 17\ar@{--}[rr]\ar[rd] && 
28 \ar@{--}[rr]\ar[rd] && 
 13 \ar[rd]  \\
 & & & 16\ar@{--}[rr]\ar[rd]\ar[ru] && 27\ar@{--}[rr]\ar[rd]\ar[ru]
 && 38\ar@{--}[rr]\ar[rd]\ar[ru] && 14\ar[rd] \\
 & & 15\ar@{--}[rr]\ar[rd]\ar[ru] && 
 26\ar@{--}[rr]\ar[rd]\ar[ru]
 && 37\ar@{--}[rr]\ar[rd]\ar[ru] &&
 48\ar@{--}[rr]\ar[rd]\ar[ru] && 
 15\ar[rd] \\
 & 14\ar@{--}[rr]\ar[rd]\ar[ru] && 25\ar@{--}[rr]\ar[rd]\ar[ru]
 && 36\ar@{--}[rr]\ar[rd]\ar[ru] &&
 47\ar@{--}[rr]\ar[rd]\ar[ru] && 58\ar@{--}[rr]\ar[rd]\ar[ru] &&
 16\ar[rd] \\
13\ar@{--}[rr]\ar[ru] && 
 24\ar@{--}[rr]\ar[ru] && 
 35\ar@{--}[rr]\ar[ru]
 && 46\ar@{--}[rr]\ar[ru] && 
 57 \ar@{--}[rr]\ar[ru] &&
 68\ar@{--}[rr]\ar[ru] && 17}
\]
Its second power has three components: one component is 
$\Gamma(2,2)$. The two other components are both quivers 
of quotients $D^b({\rm A}_3)/[1]$ of $D^b(A_3)$ by the shift. 
In particular, we have thus given a geometric construction of a quotient 
of $D^b(A_3)$ which is not an $m$-cluster category! 

The three components of the second power of $\Gamma(6,1)$ are:
\[
\small
\xymatrix@-6mm{
 & 16\ar[rdd]\ar@{--}[rr] && 38\ar[rdd]\ar@{--}[rr] && 25\ar[rdd]\ar@{--}[rr]
 && 47\ar[rdd]\ar@{--}[rr] && 16\\
 & & \\
14\ar@{--}[rr]\ar[ruu] & & 36\ar[ruu]\ar@{--}[rr] & & 58\ar[ruu]\ar@{--}[rr]
 & &  27\ar[ruu]\ar@{--}[rr] & & 14\ar[ruu]}
\]

\[
\small
\xymatrix@-8mm{
 & & 17\ar[rdd]\ar@{--}[rr] && 13\ar[rdd]\\ 
 & & \\
 & 15\ar[rdd]\ar[ruu]\ar@{--}[rr] & & 37\ar[rdd]\ar[ruu]\ar@{--}[rr] 
 & & 15\ar[rdd]\\ 
 & & \\
13\ar@{--}[rr]\ar[ruu] & & 35\ar[ruu]\ar@{--}[rr] & & 57\ar[ruu]\ar@{--}[rr]
 & &  17}
\small
\xymatrix@-8mm{ & & \\
 & & 28\ar[rdd]\ar@{--}[rr] && 24\ar[rdd]\\ 
 & & \\
 & 26\ar[rdd]\ar[ruu]\ar@{--}[rr] & & 48\ar[rdd]\ar[ruu]\ar@{--}[rr] 
 & & 26\ar[rdd]\\ 
 & & \\
24\ar@{--}[rr]\ar[ruu] & & 46\ar[ruu]\ar@{--}[rr] & & 68\ar[ruu]\ar@{--}[rr]
 & &  24}
\]
\end{example}

\begin{remark}
We have a corresponding result for type $D$, see Theorem 5.1 of~\cite{bm2}. 
In addition, in type $D$ we give an explicit description of {\em all} connected 
components appearing 
in the $m$-th power of $\Gamma_{\odot}(nm-m+1,1)$
\end{remark}





%
\section{Connections and future directions}\label{s:connection}
%
To finish we want to provide a short outlook and describe some open problems 
and possible future directions. 

\begin{itemize}
\item In recent work with Robert Marsh (\cite{bm3}) 
we provide 
a link between cluster algebra combinatorics and perfect matchings 
(for vertices and edges of a triangulation). This uses work of 
Conway-Coxeter (\cite{cc1},~\cite{cc2}), and of 
Broline-Crowe-Isaacs (\cite{bci}) on frieze patterns of positive integers. 
\item $Y$-systems can be defined in general for a pair $(G,H)$ 
of Dynkin types. 
Zamolodchikovs periodicity conjecture for general $Y$-systems 
have been proved for 
$G=A_1$ and $H=A_n$ by Frenkel-Szenes (\cite{fs}), 
by Gliozzi-Tateo (\cite{gt}) and 
for $G=A_1$, $H$ any Dynkin type by Fomin-Zelevinsky (\cite{fz-Y}) 
using cluster algebras theory. More recently, the cases 
$G=A_k$,  $H=A_n$ have been solved (\cite{sz}, \cite{vo}, independently). 
In this context, there are various open questions. First of all: can periodicity 
be proved for $G$ of arbitrary of Dynkin type using the theory 
of cluster algebras? This is not even known for $G=A_k$. 
Second: what would be a good counterpart on the side of $Y$-systems 
to the geometric model for $m$-cluster categories? 
And thirdly: In current work with Marsh we have discovered classes of infinite 
periodic systems ($G=A_1$, $H=A_{\infty}$). Does this have a translation 
to the setting of $Y$-systems? 
\item The approach to model cluster algebras with discs $(S,M)$ works 
for types $A$ and $D$ (\cite{fst06}) and for types $B$, $C$ under certain 
modifications (\cite{cfz}).  \\
{\it Open}: what can be said about the exceptional types, in particular, 
is there a way to model type $E$ using a disc with marked points?  
\item Jorgensen (\cite{jo}) has obtained 
$m$-cluster categories as quotient categories  
of cluster categories via deletion of rows ($\tau$-orbits). 
They inherit a triangulated structure. This process can be viewed as a 
reverse to our construction using the $m$-th power of a quiver. 
{\it Question}: how can we explain the triangulated structure of 
$m$-cluster categories via the $m$-th power of a quiver? 
\end{itemize}
%
%

\end{document}